\documentclass{article}
\usepackage{amssymb,amsfonts,amsmath,amsthm,amscd}
\usepackage{makeidx}
\usepackage{tikz-cd} 
\usepackage[utf8]{inputenc}
\usepackage[english]{babel}
\makeindex
\usepackage{tikz-cd}  
\usepackage[multiple]{footmisc}
\usepackage{graphicx}
\usepackage{multicol}
\usepackage{hyperref}
\usepackage{tikz}

\usepackage{hyperref}
\input xy
\xyoption{all}
\usepackage[all]{xy}

\usepackage{fancyhdr}
\usepackage{makeidx}
\makeindex 
\usepackage{vmargin} 
\setmargins{3.3cm}       
{3cm}                        
{14.5cm}                      
{23cm}                    
{0.1pt}                           
{0.5cm}                           
{0pt}                             
{1cm}                           

\numberwithin{equation}{section}

\newtheorem{teo}{Theorem}[section]

\theoremstyle{definition}

\newtheorem{lem}[teo]{Lemma}
\newtheorem{ejem}[teo]{Example}

\newcommand{\m}{{}^{-1}}

\newcommand{\ep}{\epsilon}

\newcommand{\N}{\mathbb{N}}

\title{ \textbf{ A note on adjunction spaces and \texorpdfstring{$G$-ANE's}}
}
\author{Luis Mart\'{i}nez
\footnote{ The author was supported by CONAHCYT (México)}
	 \\ 
  \small e-mail: luchomartinez9816@hotmail.com \\
	\small Departamento de Matem\'{a}ticas, Facultad de Ciencias, Universidad Nacional Aut\'{o}noma de M\'{e}xico 
}

\begin{document}

	\maketitle
	\begin{abstract} 
		We provide a different proof of the equivariant version of the Borsuk-Whitehead-Hanner Theorem in the category of proper $G$-spaces which are metrizable by a $G$-invariant metric.	\end{abstract}
 \section{Introduction} 

Given two metrizable spaces $X$ and $Y$ and a continuous map $f:A\rightarrow Y$, where $A$ is a closed subset of $X$, the Borsuk-Whitehead-Hanner Theorem (see  \cite[Chapter IV, Theorem 5.3]{Hu}, or \cite[Theorem 8.2]{Hanner1951}) states that the adjunction space $X\cup_f Y$ is ANE (absolute neighborhood extensor) for the class of metrizable spaces provided that $X$, $A$ and $Y$ are  ANE's and $X\cup_f Y$ is metrizable. More general versions of this fact were shown in \cite{Hyman1967}. In \cite[Theorem 3.11]{Antonyan2001} was shown an equivariant version of this theorem in the context of proper $G$-spaces, where $G$ is a locally compact Hausdorff group. In that proof, Hyman's idea is followed and $G$-semicanonical covers are constructed. The goal of this paper is to present a proof of the theorem without using $G$-semicanonical covers. In that sense, we present Lemma \ref{-1R} that was shown in \cite{Antonyan1990} in the context of compact groups, and from this Lemma we prove the Theorem \ref{main}. The principal tool of the proof is to construct equivariant deformations of the adjunction space.

 \section{The notions}\label{pre}
Throughout the paper the letter $G$ will always denote a locally compact Hausdorff group and $\mathbb{N}=\{1,2,\dots\}$ denotes the natural numbers. 

By a topological transformation group we mean a triple $(G,X,\theta)$, where $\theta: G\times X\rightarrow X$ is a continuous action, that is, a continuous map such that $\theta(g,\theta(h,x))=\theta(gh,x)$ and $\theta(e,x)=x$ for any $g, h\in G$ and any $x\in X$, where $e$ is the unity of $G$. A space $X$ together with
a fixed action $\theta$ of the group $G$ is called a $G$-space.  It is usual to write $gx$ instead of $\theta(g,x)$ for the image $\theta(g,x)$ of a pair $(g,x)\in G\times X$. In a similar way, we write $G(S)=\theta(G\times S)$ for each $S\subseteq X$. In particular,
$G(x)$ denotes the $G$-orbit of $x\in X$. For any $x\in X$ we denote $G_x=\{g\in G: gx=x\}$ the isotropy subgroup (or stabilizer) of $x\in X$.

The orbit space $X/G=\{G(x): x\in X\}$ is endowed with the quotient topology determined by the orbit map $p :X\rightarrow X/G$, $x\mapsto G(x)$, $x\in X$. If $S\subseteq X$ satisfies that $S=p^{-1}(p(S))$ (that is, $S=G(S)$), then $S$ is called a $G$-subset (or an invariant subset) of $X$.

If $Y$ is another $G$-space, a continuous map $f: X\rightarrow Y$ is called a G-map (or an equivariant map) if $f(gx) = gf(x)$ for every $g\in G$ and $x\in X$. If $G$ acts trivially on $Y$ then $f$ is called invariant map.

Let $X$ and $Y$ be $G$-spaces. Two $G$-maps $p, q: X\rightarrow Y$ are called $G$-homotopic provided that there exists a continuous map $F: X\times I\rightarrow Y$ such that $F_t: X\rightarrow Y$, $x\mapsto F(x,t)$, is a $G$-map for any $t\in I$, $F_0=p$ and $F_1=q$. In this case $F$ is called a $G$-homotopy. Take $\mathcal{V}$ an open cover of $Y$.  $F$ is said to be \textit{limited} by $\mathcal{V}$ (or a $(\mathcal{V},G)$-homotopy) provided that for any $x\in X$ there exists $V\in \mathcal{V}$ such that $F(\{x\}\times I)\subseteq V$. In this case $p$ and $q$ are called $(\mathcal{V},G)$-homotopic.  Let $\ep>0$ be a real number and suppose that $Y$ is metrizable. $F$ is called $(\ep,G)$-homotopy provided that $F(\{x\}\times I)$ has diameter less than $\ep$, for each $x\in X$. 

Take $X$ a $G$-space. A $G$-homotopy $F:X\times I\rightarrow X$ such that $F_0$ is the identity map on $X$ is called a $G$-deformation of $X$. If $\mathcal{V}$ is an open cover of $X$ and $F$ is $(\mathcal{V},G)$-homotopy, then $F$ is called $(\mathcal{V},G)$-deformation. If $X$ is metrizable and $F$ is $(\varepsilon,G)$-homotopy, then $F$ is called $(\varepsilon,G)$-deformation.  A sequence of $G$-deformations $(F^n)_{n\in \N}$ of $X$ is said to converge to the identity map, provided that for every $x\in X$ and every neighborhood $V$ of $x$ in $X$, there exists a neighborhood $W$ of $x$ in $X$ and $k\in \N$ so that $F^n(W\times I)\subseteq V$, for every $n\geq k$.
\begin{ejem}\label{ejeconverid}
   Let $(X,d)$ be a metric $G$-space. Suppose that for each $n\in\N$ there exists a $(\frac{1}{n},G)$-deformation $F^n$ of $X$.
   We are going to check that $(F^n)_{n\in \N}$ converges to the identity. Take $x\in X$ and let $V$ be a neighborhood of $x$ in $ X$. There exists $n\in\N$ such that $B_d(x,\frac{1}{n})\subseteq V$. Observe that $F^k(B_d(x,\frac{1}{2n+1})\times I)\subseteq V$, for each $k\geq 2n+1$. Indeed, take $k\geq 2n+1$, $y\in B_d(x,\frac{1}{2n+1})$ and $t\in I$. Then:
\begin{equation*}
        \begin{aligned}
      d(F^k(y,t), x) & \leq d(F^k(y,t),F^k(y,0))+d(F^k(y,0),x) \\
    &\leq\mbox{diam}(F^k(\{y\}\times I))+d(y,x)\\
    & 
    <\frac{1}{k}+\frac{1}{2n+1}\\
    & <\frac{1}{n},
        \end{aligned}
    \end{equation*}
then $F^k(y,t)\in V$ and we conclude that $F^k(B_d(x,\frac{1}{2n+1})\times I)\subseteq V$. This shows that $F$ converges to the identity.
\end{ejem}
Let $X$ be a $G$-space. A subset $V\subseteq X$ is called small, if for every point of $X$ there is a neighborhood $U$ in $X$ with the property that the set $\langle U,V\rangle=\{g\in G: gU\cap V\neq \emptyset\}$ has compact closure in $G$. Recall that  $X$ is called proper (in the sense of Palais \cite[Definition 1.2.2]{Palais61}), if $X$ has an open cover consisting of small sets. In this case, each stabilizer is compact and each orbit
is closed \cite[Proposition 1.1.4]{Palais61}.

In the present paper we are interested in the class $G$-$\mathcal{M}$ of all metrizable proper $G$-spaces $X$ that admit a $G$-invariant
metric. Remember that a compatible metric $d$ on a proper $G$-space $X$ is called $G$-invariant provided that $d(gx,gy)=d(x,y)$ for any $x, y\in X$ and $g\in G$. In this case the orbit space is metrizable. Indeed, the function
\begin{center}
    $\tilde{d}(G(x),G(y))=\mbox{inf}\{d(\tilde{x},\tilde{y}): \tilde{x}\in G(x), \tilde{y}\in G(y)\}$,
\end{center}
defines a compatible metric with the quotient topology of $X/G$ \cite[Theorem 4.3.4]{Palais61}.

Let $Y$ be a $G$-space. A closed $G$-subset $A$ of $Y$ is called $G$-neighborhood retract (resp. $G$-retract) of $Y$ if there exists a $G$-neighborhood $U$ of $A$ (resp. $U=Y$) and a $G$-retraction $r: U\rightarrow A$ (resp. $r: Y\rightarrow A$). 

A $G$-space $X$ is called \textit{absolute neighborhood extensor} (resp. \textit{absolute extensor}) of $G$-$\mathcal{M}$ provided that, for every $Y\in G$-$\mathcal{M}$ and every closed $G$-subset $A$ of $Y$, every $G$-map $f: A\rightarrow X$ can be extended to a $G$-map $F: U\rightarrow X$ on some $G$-neighborhood $U$ of $A$ in $Y$ (resp. on all of $Y$). In this case we write $X\in G$-A(N)E.

A $G$-space $X\in G$-$\mathcal{M}$ is called \textit{absolute neighborhood retract} (resp. \textit{absolute retract}) of $G$-$\mathcal{M}$  provided that, for every $Y\in G$-$\mathcal{M}$ and every closed $G$-embedding $\iota: X\rightarrow Y$, the image $\iota(X)$ is $G$-retract of some  $G$-neighborhood in $Y$ (resp. of all of $Y$). In this case we write $X\in G$-A(N)R.

Given $X\in G$-$\mathcal{M}$, it is known by \cite[Corollary 6.3]{Antonyan2014} that $X\in G$-A(N)E if and only if $X\in G$-A(N)R.

In the sequel we will need the following results:

\begin{teo}\label{embedingAE}
    \cite[Theorem 6.1]{Antonyan2014} For every $X\in G$-$\mathcal{M}$ there exists $L\in G$-$\mathcal{M}$ such that $L\in G$-AE and $X$ can be considered as an invariant closed subset of $L$.
\end{teo}
\begin{lem}\label{twopro}
    Let $X\in G$-$\mathcal{M}$. Then:
    \begin{enumerate}
        \item [(i)] If $A$ and $B$ are disjoint closed $G$-subsets of $X$, then there exists an invariant map $F: X\rightarrow I$ such that $F(A)=\{0\}$ and $F(B)=\{1\}$. In particular, there are disjoint open $G$-subsets $U$ and $V$ of $X$ such that $A\subseteq U$ and $B\subseteq V$.
        \item [(ii)] Let $A$ be a  closed $G$-subset of $X$ and let $O$ be an open $G$-subset of $X$. If $A\subseteq O$, then there exists $W$ an open $G$-subset of $X$ such that $A\subseteq W\subseteq \overline{W}\subseteq O$.
    \end{enumerate}
\end{lem}
\begin{proof}
    (i) Let $p: X\rightarrow X/G$ be the orbit map. Observe that $A=p\m( p (A))$ and $B=p\m( p(B))$, then $p(A)$ and $p(B)$ are disjoint closed subsets of $X/G$. Since $X/G$ is normal, then there exists an invariant map $f: X/G\rightarrow I$ such that $f(p(A))=\{0\}$ and $f(p(B))=\{1\}$.  Notice that $F=f\circ p$ satisfies the desired condition. 

    (ii) Put $B=X\setminus O$. Since $A$ and $B$ are disjoint closed $G$-subsets of $X$, by (i) there are disjoint open $G$-subsets $W$ and $V$ of $X$ such that $A\subseteq W$ and $B\subseteq V$. Then $A\subseteq W\subseteq\overline{W}\subseteq X\setminus V\subseteq O$. 
\end{proof}
\begin{lem}\label{dowker}

Let $X, Y\in G$-$\mathcal{M}$ and let $A$ be a closed $G$-subset of $X$. Let $f: T\rightarrow Y$ be a $G$-map, where $T$ is the  closed $G$-subset $T=(X\times \{0\})\cup (A\times I)$ of $X\times I$. If $f$ can be extended to a $G$-map $\tilde{f}: (X\times \{0\})\cup U\rightarrow Y$ for some $G$-neighborhood $U$ of $A\times I$ in $X\times I$, then $f$ extends to a $G$-map $F: X\times I\rightarrow Y$. 
\end{lem}

\begin{proof}
 Since $I$ is compact, then there exists  a neighborhood $\tilde{V}$ of $A$ in $X$ such that $\tilde{V}\times I\subseteq U$. Since $U$ is invariant, then $V\times I\subseteq U$, where $V=G(\tilde{V})$. Notice that $V$ is a $G$-neighborhood of $A$ in $X$. By Lemma \ref{twopro} there exists an invariant map $e: X\rightarrow I$ such that $e(A)=\{1\}$ and $e(X\setminus V)=\{0\}$.

 Consider $F: X\times I\rightarrow I$, given by:
 \begin{center}
     $F(x,t)=\tilde{f}(x,e(x)t)$, for each $(x,t)\in X\times I$,
 \end{center}
 which is a $G$-map and extends to $f$.
\end{proof}

\begin{lem}\label{lemrede} \cite[Lemma 3.10]{Antonyan2001}
    Let $X\in G$-ANE $\cap\ G$-$\mathcal{M}$. If $A$ is a closed $G$-subset of $X$ and $A\in G$-ANR, then for any open cover $\mathcal{V}$ of $X$ there exists a $(\mathcal{V},G)$-homotopy $H:X\times I \rightarrow X$ such that:
    \begin{enumerate}
        \item [(i)] $H_0=\mbox{id}_X$.
        \item [(ii)] $H(a,t)=a$, for any $a\in A$ and any $t\in I$.
        \item [(iii)] There exists a $G$-neighborhood $V$ of $A$ in $X$ such that $H_1(V)=A$.
    \end{enumerate}
\end{lem}
\begin{lem}\label{restrong}
    Let $X\in G$-AE and let $A$ be a closed $G$-subset of $X$. Then $A\in G$-AE if and only if $A$ is a  strong $G$-deformation retract of $X$. 
\end{lem}
\begin{proof}
    $\Rightarrow)$ Consider the $G$-map $f: (X\times\{0\})\cup (A\times I)\rightarrow A$ given by:
    \begin{center}
  $f(x,t)=\left\{
	\begin{array}{ll}
		 x, & \mbox{if } x\in X,\ t=0 \\
  
		x, & \mbox{if } x\in A,\ t\in I
	\end{array}
\right. 
        $
   
\end{center}
Observe $X\times I\in G$-$\mathcal{M}$ and $A\in G$-AE, then $f$ can be extended to a $G$-map $F: X\times I\rightarrow A$. Notice that $F_1$ is a $G$-retraction.

 $\Leftarrow)$ Since $A$ is $G$-retract of a $G$-AE, then $A\in G$-AE.
\end{proof}

 \section{Adjunction spaces and $G$-ANE's}
\begin{lem}\label{-1R}
    Let $X\in G$-$\mathcal{M}$ and let $L\in G$-$\mathcal{M}\cap G$-AE such that $X$ is a closed $G$-subset of $L$. Then the following statements are equivalent:
    \begin{enumerate}
        \item [(i)] $X\in G$-ANE.
        \item [(ii)] For each open cover $\mathcal{V}$ of $X$ there exists a $(\mathcal{V},G)-$deformation $F:X\times I\rightarrow X$ of $X$ such that $F_1$ is equivariantly extendable to some $G$-neighborhood of $X$ in $L$.
        \item [(iii)] For some metric on $X$ and for any $\ep>0$, there exists a $(\ep,G)$-deformation $F:X\times I\rightarrow X$ such that $F_1$ is equivariantly extendable to some $G$-neighborhood of $X$ in $L$.
        \item [(iv)] There exists a sequence of $G$-deformations $(F^n)_{n\in \N}$ of $X$ converging to the identity map such that $F^n_1$ extends to some $G$-neighborhood of $X$ in $L$, for any $n\in \N$.  
        
    \end{enumerate}
    \end{lem}
   \begin{proof}
(i) $\Rightarrow$ (ii)  Consider $\mathcal{V}$ an open cover of $X$ and define $F:X\times I\rightarrow X$ given by $F(x,t)=x$, $x\in X$, $t\in I$. Put $F^n=F$ for any $n\in \N$ and observe that $(F^n)_{n\in \N}$ satisfies (ii) because $X\in G$-ANE.

(ii)$\Rightarrow$ (iii) Let $d$ be a compatible metric on $X$. Take $\ep>0$ and put $\mathcal{V}_\ep=\{B_d(x,\frac{\ep}{2}): x\in X\}$ which is an open cover of $X$. By (i) there exists a $(\mathcal{V},G)$-deformation $F$ of $X$ such that $F_1$ is equivariantly extendable to some $G$-neighborhood of $X$ in $L$. Notice that $F$ is a $(\ep,G)$-deformation.

(iii)$\Rightarrow$(iv) For each $n\in \mathbb{N}$ let $F^n$ be a $(\frac{1}{n},G)$-deformation of $X$ such that $F^n_1$ can be extended to a $G$-neighborhood of $X$ in $L$. It was shown in Example \ref{ejeconverid} that $(F^n)_{n\in \mathbb{N}}$ converges to the identity.

(iv) $\Rightarrow$ (i) Take $U$ a $G$-neighborhood of $X$ in $L$ and let $f: U\rightarrow X$ be a $G$-extension of $F^1_1.$ Put $t_n=1-2^{-(n-1)}$ for any $n\in \N$. We want to construct a $G$-map $F: U\times [0,1)\rightarrow X$ such that $F(x,t_n)=F^n_1(x)$, for each $x\in X$ and each $n\in \N$.

First define $F$ on $U\times \{0\}$ by $F(x,0)=f(x)$, $x\in U$. Since $t_1=0$, suppose that $F$ has been defined on $U\times [t_1,t_n]$ for some $n\geq 1$. Now we shall construct $F$ on $U\times [t_n,t_{n+1}]$. To do that, consider the closed $G$-subset $T= (X\times I)\cup (U\times \{0\})$ of $U\times I$, and observe that $k: T\rightarrow X$ defined by:
\begin{center}
  $k(x,t)=\left\{
	\begin{array}{ll}
		 F(x,t_n), & \mbox{if } x\in U,\ t=0 \\
  
		F^n(x,1-t), & \mbox{if } x\in X,\ t\in I
	\end{array}
 \right. 
       $
    
 \end{center}
is a $G$-map. Indeed, take $g\in G$ and $(x,t)\in T$. If $x\in U$ and $t=0$, then:
\begin{center}
    $k(g(x,t))=k(gx,0)=F(gx,t_n)=gF(x,t_n)=gk(x,t)$,
\end{center}
in a similar way, if $x\in X$, $t\in I$, then:
\begin{center}
    $ k(g(x,t))=k(gx,t)=F^n(gx,1-t)=gF^n(x,1-t)=gk(x,t)$. 
\end{center}
This shows that $k$ is a $G$-map. 

Since $F_1^{n+1}: X\rightarrow X$ can be extended to a $G$-neighborhood of $X$ in $L$ and $L\in G$-AE, then $F_1^{n+1}\circ k: T\rightarrow X$ extends to a $G$-neighborhood of $T$ in $U\times I$. Follows by Lemma \ref{dowker} that there exists a $G$-extension  $K: U\times I\rightarrow X$ of $F_1^{n+1}\circ k$.

Define $F$ on $U\times [t_n,t_{n+1}]$ by:

\begin{center}
  $F(x,(1-t)t_n+tt_{n+1})=\left\{
	\begin{array}{ll}
		 F_{2t}^{n+1}(F(x,t_n)), & \mbox{if } 0\leq t\leq \frac{1}{2} \\
  
		K_{2t-1}(x), & \mbox{if } \frac{1}{2}\leq t\leq 1
	\end{array}
\right. 
        $
    
\end{center}
Notice that $F$ is a $G$-map. Indeed, take $g\in G$ and $(x,s)\in U\times [t_n,t_{n+1}]$. Then $s=(1-t)t_n+tt_{n+1}$ for some $0\leq t\leq 1 $. If $0\leq t\leq \frac{1}{2}$, then:
\begin{center}
    $F(gx,s)=F(gx,(1-t)t_n+tt_{n+1})=F_{2t}^{n+1}(F(gx,t_n))=gF_{2t}^{n+1}(F(x,t_n))=gF(x,s)$.
\end{center}
 On the other hand, if $\frac{1}{2}\leq t\leq 1$, then:
 \begin{center}
   $F(gx,s)=F(gx,(1-t)t_n+tt_{n+1})=K(gx,2t-1)=gK(x,2t-1)=gF(x,s)$,  
 \end{center}
 then we are done. In this way we have constructed the $G$-map $F:U\times [0,1)\rightarrow X$.

 Now define the $G$-map $H:X\times I\rightarrow X$ by $H\restriction_{X\times [0,1)}=F\restriction_{X\times [0,1)}$ and $H(x,1)=x$ for every $x\in X$. To show that $H$ is continuous it only remains to verify the continuity on $X\times \{1\}$. Take $x_0\in X$ and let $V$ be a neighborhood of $H(x_0,1)=x_0$ in $X$. Since $(F^n)_{n\in \N}$ converges to the identity, there are neighborhoods $W$ and $W_1$ of $x_0$ in $X$ and $k, k_1\in \mathbb{N}$ such that $k\geq k_1$ and

\begin{center}
    $F^{n+1}(W_1\times I)\subseteq V$, for each $n\geq k_1$, and $F^n(W\times I)\subseteq W_1$, for each $n\geq k$.
\end{center}
We shall check that $H(W\times [t_k,1])\subseteq V$. Take $x\in W$, $s\in [t_k,1]$, $n\geq k$ and $t\in [0,1]$ such that $s\in [t_n,t_{n+1}]$ and $s=(1-t)t_n+tt_{n+1}$. Then $\{F^n_1(x), F^n_{2-2t}(x)\}\subseteq W_1$ and $\{F^{n+1}_{2t}(F^n_1(x)), F^{n+1}_1(F^n_{2-2t}(x))\}\subseteq V$. But
\begin{center}
  $H(x,s)=\left\{
	\begin{array}{ll}
		 F_{2t}^{n+1}(F^n_1(x)), & \mbox{if } 0\leq t\leq \frac{1}{2} \\
  
		F^{n+1}_1(F^n_{2-2t}(x)), & \mbox{if } \frac{1}{2}\leq t\leq 1
	\end{array}
\right. 
        $
    
\end{center}
then $H(x,s)\in V$.

On $U\times [0,1)$ consider the metric $\rho$ given by $\rho((u_1,t_1),(u_2,t_2))=d(u_1,u_2)+|t_1-t_2|$, where $d$ is an invariant metric on $L$. Consider the $G$-neighborhood of $X\times [0,1)$ given by:
\begin{center}
    $O=\left\lbrace 
    (u,t)\in U\times [0,1): \begin{array}{c} \exists (x,s)\in X\times [0,1),  \\
    0\leq \rho((u,t),(x,s)),d(F(u,t),F(x,s)) <1-t
    \end{array}\right\rbrace
 $.
\end{center}
We shall check that $O$ is invariant. Take $(u,t)\in O$ and $g\in G$. There exists $(x,s)\in X\times [0,1)$ such that $0\leq \rho((u,t),(x,s)),\rho(F(u,t),F(x,s)) <1-t$. Since $X$ is invariant, then $(gx,s)\in X\times [0.1)$. Moreover,
\begin{center}
    $\rho((gu,t),(gx,s))=d(gu,gx)+|t-s|=d(u,x)+|t-s|=\rho((u,t),(x,s))<1-t$,
\end{center}

\noindent in a similar way $d(F(gu,t),F(gx,s))<1-t$, then $(gu,t)\in O$ and $O$ is invariant.

Notice that $X\times [t_n,t_{n+1}]\subseteq O$, for every $n\geq 1$.  For $n=1$, since $X\times [t_1,t_2]\subseteq O$, then there exists a neighborhood $\tilde{U}_1$ of $X$ in $U$ such that $\tilde{U}_1\times [t_1,t_2]\subseteq O$. But $O$ is invariant, then $U_1\times [t_1,t_2]\subseteq O$, where $U_1=G(\tilde{U}_1)\subseteq U$ is a $G$-neighborhood of $X$ in $U$.  In a similar way, for $n=2$ there exists a $G$-neighborhood $\tilde{U}_2$ of $X$ in $U$ such that $\tilde{U}_2\times [t_2,t_3]\subseteq O$.  It follows by Lemma \ref{twopro} that there exists a $G$-neighborhood $U_2$ of $X$ in $U$ such that $U_2\subseteq \overline{U}_2\subseteq U_1\cap \tilde{U}_2$, then $U_2\times [t_2,t_3]\subseteq O$. In this way, inductively there is a $G$-neighborhood $U_n$ of $X$ in $U$ such that $U_n\times [t_n,t_{n+1}]\subseteq O$ and $\overline{U_{n+1}}\subseteq U_n$, $n\in \mathbb{N}$.

For each $n\geq 1$ let $e_n: U\rightarrow I$ be an invariant map such that:
\begin{center}
  $e_n(u)=\left\{
	\begin{array}{ll}
		 0, & \mbox{if } u\in U\setminus U_n \\
 
		1, & \mbox{if } u\in \overline{U_{n+1}}
	\end{array}
\right. 
       $
    
\end{center}
then $e:U\rightarrow I$, $u\mapsto \sum\limits_{n=1}^\infty\frac{e_n(u)}{2^n}$, $u\in U$, is an invariant map. 

Take $u\in U_1\setminus X$. We shall check that $(u,e(u))\in O$. Since $X$ is closed in $L$, then $X=\cap_{n=1}^\infty U_n$ and there exists $m\geq 1$ such that $u\in U_m\setminus U_{m+1}$. Then:

\begin{center}
  $e_n(u)=\left\{
	\begin{array}{ll}
		 1, & \mbox{if } n<m \\
  
		0, & \mbox{if } n>m\
	\end{array}
\right. 
        $
    
\end{center}
and
\begin{center}
    $e(u)=\frac{1}{2}+\cdots+\frac{1}{2^{m-1}}+\frac{1}{2^m}e_m(u)=1-\frac{1}{2^{m-1}}+\frac{1}{2^m}e_m(u)=t_m+\frac{1}{2^m}e(u)\in [t_m,t_{m+1}]$,
\end{center}
hence $(u,e(u))\in U_m\times [t_m,t_{m+1}]\subseteq O$.

Consider $r: U\rightarrow X$, given by:
\begin{center}
  $r(u)=\left\{
	\begin{array}{ll}
		 F(u,e(u)), & \mbox{if } u\in U\setminus X \\
  
		u, & \mbox{if } u\in X\
	\end{array}
\right. 
        $
   
\end{center}

To prove that $r$ is continuous, take $x\in \mbox{Fr}_U(X)$ and let  $(u_n)_{n\in \N}$ be a sequence in $U\setminus X$ such that $\mbox{lim}(u_n)=x$. We shall check that $\lim(r(u_n))=r(x)=x$. Since $x\in U_1$, we can suppose that $u_n\in U_1$, for each $n\in \mathbb{N}$. Take $n\in \N$. We know that $(u_n,e(u_n))\in O$, then there exists $(x_n,s_n)\in X\times [0,1)$ such that:
\begin{center}
    $\rho((u_n,e(u_n)),(x_n,s_n))<1-e(u_n)$ and $d(F(u_n,e(u_n)),F(x_n,s_n))<1-e(u_n)$.
\end{center}

But $\lim(u_n)= x$, then $\lim(e(u_n))=e(x)=1$ and $\lim(u_n,e(u_n))=(x,1)$. This shows that:

\begin{center}
   $\lim(\rho((u_n,e(u_n)),(x_n,s_n)))=\mbox{lim}(d(F(u_n,e(u_n)),F(x_n,s_n)))=0$,
\end{center}
then $\lim(x_n,s_n)=(x,1)$ and  $\lim(H(x_n,s_n))=H(x,1)=x$. But $H(x_n,s_n)=F(x_n,s_n)$ for each $n\in \N$, then: 
\begin{center}
    $\lim(r(u_n))=\lim(F(u_n,e(u_n)))=\lim(F(x_n,s_n))=\lim(H(x_n,s_n))=x$,
\end{center}
hence $r$ is continuous. Observe that $r$ is a $G$-map, then $X$ is a $G$-neighborhood retract of $U$. We conclude that $X\in G$-ANE.

\end{proof}

Consider $(G,X,\alpha)$ and $(G,Y,\beta)$ topological transformation groups.  Let $A$ be a closed $G$-subset of $X$ and let $f: A\rightarrow Y$ be a $G$-map. Remember that the \textit{adjunction space} obtained by adjoining $X$ to $Y$ by means of $f$, denoted by $X\cup_f Y$, is the quotient space of $(X\sqcup Y)/\sim$, where $\sim$ is the equivalence relation generated by:
    \begin{center}
        $\{(a,f(a))\in (X\sqcup Y)\times (X\sqcup Y): a\in A\}$.
    \end{center}
    There is a canonical action of $G$ on $X\cup_f Y$, defined as follows: first consider the action $\alpha\sqcup\beta: G\times (X\sqcup Y)\rightarrow X\sqcup Y$, given by:
    \begin{center}
        $(\alpha\sqcup\beta) (g,z)=\left\{
	\begin{array}{ll}
		 \alpha(g,z), & \mbox{if } z\in X\\
  
		\beta(g,z), & \mbox{if } z\in Y
	\end{array}
\right. 
        $
    \end{center}
    and let $p: X\sqcup Y\rightarrow X\cup_f Y$ the quotient map. Since $G$ is locally compact, then $\mbox{id}_G\times p$ is an identification and by the universal property of a quotient, there exists a continuous map $\theta: G\times X\cup_fY\rightarrow X\cup_f Y$ such that the following diagram is commutative: 
 
    \begin{center}
         $\xymatrix{G\times (X\sqcup Y)\ar[d]_-{\mbox{id}_G\times p}\ar[r]^-{\alpha\sqcup \beta}& X\sqcup Y\ar[d]^-{p}\\
         G\times (X\cup_f Y) \ar[r]_-{\theta}&  X\cup_f Y
    }$
    \end{center}
    In the next result $X\cup_f Y$ is considered  as a $G$-space with the action $\theta$.

\begin{teo}\label{main}
    Let $X, Y\in G$-$\mathcal{M}$, let $A$ be a closed $G$-subset of $X$ and let $f: A\rightarrow Y$ be a $G$-map. If $X$, $A$ and $Y$ are $G$-ANE's and $X\cup_f Y\in G$-$\mathcal{M}$, then $X\cup_f Y$ is $G$-ANE.
\end{teo}
\begin{proof}
    Let $\mathcal{V}$ be an open cover of $X\cup_f Y$. Notice that $\mathcal{U}=p\m(\mathcal{V})=\{p\m(V): V\in \mathcal{V}\}$ is an open cover of $X\sqcup Y$. Since $\mathcal{W}=\{X\cap U: U\in\mathcal{U}\}$ is an open cover of $X$, by Lemma \ref{lemrede} there exists a $(\mathcal{W},G)$-homotopy $T: X\times I\rightarrow X$ such that $T_0=\mbox{id}_X$, $T(a,t)=a$ for each $a\in A$ and $t\in I$, and there exists a $G$-neighborhood $V$ of $A$ in $X$ so that $T_1(V)=A$.

    Consider the $(\mathcal{U},G)$-homotopy $K: (X\sqcup Y)\times I\rightarrow X\sqcup Y$ given by:
    \begin{center}
  $K(z,t)=\left\{
	\begin{array}{ll}
		 T(z,t), & \mbox{if } z\in X \\
  
		z, & \mbox{if } z\in Y
	\end{array}
\right. 
        $
    
\end{center}
There is a continuous map $H: (X\cup_fY)\times I\rightarrow X\cup_fY$ such that the following diagram is commutative:

\begin{center}
         $\xymatrix{ (X\sqcup Y)\times I\ar[d]_-{p\times \mbox{id}_I}\ar[r]^-{K}& X\sqcup Y\ar[d]^-{p}\\
          (X\cup_f Y)\times I \ar[r]_-{H}&  X\cup_f Y
    }$
    \end{center}

    Observe that $H$ is a $(\mathcal{V},G)$-deformation of $X\cup_fY$. Indeed, take $z\in X\cup_f Y$, $t\in I$, $g\in G$ and $w\in X\sqcup Y$ such that $p(w)=z$. Since $p(gw)=gz$, then:
    \begin{center}
        $H(gz,t)=p(K(gw,t))=p(gK(w,t))=gp(K(w,t))=gH(z,t)$,
    \end{center}
    hence $H$ is a $G$-map. Now take $U\in \mathcal{U}$ and $V\in \mathcal{V}$ such that $p(U)\subseteq V$ and $K(\{w\}\times I)\subseteq U$. Then $H(\{z\}\times I)=p(K(\{w\}\times I))\subseteq V$, and this shows that $H$ is limited by $\mathcal{V}.$

    By Theorem \ref{embedingAE}  there exists $L\in G$-AE such that $X\cup_fY$ can be considered as a closed $G$-subset of $L$. We shall check that $H_1$ can be extended to a $G$-neighborhood of $X\cup_f Y$ in $L$.

    Since $p\restriction_Y: Y\rightarrow p(Y)$ is a closed $G$-equivalence, we identify $Y=p(Y)$ as a $G$-space of $X\cup_f Y$. Observe that $F=p(X\setminus V)$ is a closed $G$-subset of $L$ and $F\cap Y=\emptyset$. By Lemma \ref{twopro} there are $R$ and $S$ open $G$-subsets of $L$ so that $F\subseteq R$, $Y\subseteq S$ and $R\cap S=\emptyset$. Consider the following $G$-subsets of $L$: $B=L\setminus (R\cup S)$, $C=B\cap (X\cup_f Y)\subseteq p(V\setminus A)$, $D=R\cap(X\cup_fY)\subseteq p(X\setminus A)$ and $E=S\cap (X\cup_f Y)\subseteq p(V\cup Y).$

    Note that $\phi: C\rightarrow A$ given by $\phi(c)=K_1(p\m(c))$ is a $G$-map. Indeed, take $c\in C$, $g\in G$ and $v\in V\setminus A$ such that $p(v)=c$. Then $p(gv)=gc$ and:
    \begin{center}
        $\phi(gc)=K_1(gv)=gK_1(v)=gK_1(p\m(c))=g\phi(c)$.
    \end{center}
    Since $A\in G$-ANE and $C$ is a closed $G$-subset of $B$, then there exists a $G$-extension $\tilde{\phi}: W_0\rightarrow A$ of $\phi$, where $W_0$ is a $G$-neighborhood of $C$ in $B$.

    Since $W_0$ and $C\cup D$ are closed $G$-subsets of $W_0\cup D$ and $W_0\cap (C\cup D)= C$, then the following map is equivariant: $\tilde{\varphi}_1: W_0\cup D\rightarrow X$, given by:

        \begin{center}
  $\tilde{\varphi}_1(s)=\left\{
	\begin{array}{ll}
		 \tilde{\phi}(s), & \mbox{if } s\in W_0 \\
  
		K_1(p\m(s)), & \mbox{if } s\in C\cup D
	\end{array}
\right. 
        $
    
\end{center}
Moreover, since $X$ is a $G$-ANE and $W_0\cup D$ is a closed $G$-subset of $W_0\cup R$ , then there exists a $G$-neighborhood $W_1$ of $W_0\cup D$ in $W_0\cup R$ and a $G$-extension $\varphi_1: W_1\rightarrow X$  of $\tilde{\varphi}_1$.

Notice that $W_0\cap (C\cup E)=C$ and $p(\tilde{\phi}(s))=p(\phi(s))=pK_1p\m(s)=H_1(s)$, for any $s\in C$. Moreover, $W_0$ and $C\cup E$ are closed $G$-subsets of $W_0\cup E$. Since $\tilde{\phi}$ and $p$ are $G$-maps, then $\tilde{\varphi}_2: W_0\cup E\rightarrow Y$ given by

        \begin{center}
  $\tilde{\varphi}_2(s)=\left\{
	\begin{array}{ll}
		 p(\tilde{\phi}(s)), & \mbox{if } s\in W_0 \\
  
		H_1(s), & \mbox{if } s\in C\cup E
	\end{array}
\right. 
        $
    
\end{center}
is a $G$-map. But $Y\in G$-ANE and $W_0\cup E$ is a closed $G$-subset of $W_0\cup S$, then there exists $\varphi_2: W_2\rightarrow Y$ a $G$-extension of $\tilde{\varphi}_2$, where $W_2$ is a $G$-neighborhood of $W_0\cup E$ in $W_0\cup S$.

Take $s\in W_1\cap W_2=W_0$. Notice that $p(\varphi_1(s))=p(\tilde{\phi}(s))=\varphi_2(s)$. Then the map $\varphi: W\rightarrow Z$ given by:

        \begin{center}
  $\varphi(s)=\left\{
	\begin{array}{ll}
		 p(\varphi_1(s)), & \mbox{if } s\in W_1 \\
  
		p(\varphi_2(s)), & \mbox{if } s\in W_2
	\end{array}
\right. 
        $
    
\end{center}
is well defined, where $W=W_1\cup W_2$. Since $W_1=W\setminus S$ and $W_2=W\setminus R$, then $\varphi$ is continuous. Moreover, $\varphi$ is a $G$-map because $p\circ \varphi_1$ and $p\circ\varphi_2$ are $G$-maps.

Notice that $W$ is a $G$-neighborhood of $X\cup_f Y$ in $L$. To prove that $W$ is open, observe that $W_1$ is open in $B\cup R$ and $B\cup R$ is closed in $L$, then $(B\cup R)\setminus W_1$ is closed in $L$. Similarly $(B\cup R)\setminus W_2$ is closed in $L$, then $L\setminus W=[(B\cup R)\setminus W_1] \cup [(B\cup S)\setminus W_2]$ is closed in $L$ and $W$ is open.

It only remains to verify that $\varphi$ is an extension of $H_1$. Take $s\in X\cup_f Y$. If $s\in C\cup D$, then $s\in W_1$ and:
\begin{center}
    $\varphi(s)=p(\varphi_1(s))=p(\phi(s))=p(K_1(p\m(s)))=H_1(s)$.
\end{center}
If $s\in E,$ then $s\in W_2$ and: 
\begin{center}
    $\varphi(s)=p(\varphi_2(s))=\varphi_2(s)=H_1(s)$.
\end{center}
Then by Lemma \ref{-1R} we conclude that $X\cup_fY$ is $G$-ANE.

\end{proof}

  
    


	\end{document}